\documentclass{amsart}

\newtheorem{theorem}{Theorem}

\newtheorem{lemma}{Lemma}

\theoremstyle{definition}
\newtheorem{example}{Example}

\newcommand{\ol}{\overline}

\title{On cyclic numbers and an extension of Midy's theorem}
\author{Juan B. Gil \and Michael D. Weiner}

\begin{document}

\maketitle

In this note we consider fractions of the form $\frac{1}{m}$ and their 
floating-point representation in various arithmetic bases. For instance,
what is $\frac17$ in base 2005? And, what about $\frac14$? We give a
simple algorithm to answer these questions.

In addition, we discuss an extension of Midy's theorem whose proof relies 
on elementary modular arithmetic. 

\section{Cyclic numbers and change of base}
Let us start with the simple and commonly used example $p=7$.  
The number $\frac17=0.\ol{142857}$ has a couple of fascinating 
properties that can be used to delight friends, even if they are
familiarized with the mysteries of math. With the period 142857 we can 
associate the ``key'' $\langle 132645\rangle$ which in this particular
case represents the order of the digits 
\begin{displaymath} 
\begin{array}{cccccc}
1& 4& 2& 8& 5& 7\\
\uparrow&\uparrow&\uparrow&\uparrow&\uparrow&\uparrow \\
1& 3& 2& 6& 4& 5
\end{array}
\end{displaymath}
so it indicates that 1 is the first digit, 4 is the third, 2 is the 
second, and so on. More precisely, in mathematical terms, the key consists 
of the residues  $10^i$ mod 7 for $i=0,1,\dots,5$. That is, 
\begin{alignat*}{3}
  10^0 &\equiv \mathbf{1} \mod 7,\quad
  & 10^1 &\equiv \mathbf{3} \mod 7,\quad
  & 10^2 &\equiv \mathbf{2} \mod 7,\\
  10^3 &\equiv \mathbf{6} \mod 7,\quad
  & 10^4 &\equiv \mathbf{4} \mod 7,\quad
  & 10^5 &\equiv \mathbf{5} \mod 7.
\end{alignat*}
If we have the period and the key, we can perform some nice computations.
First of all, we can immediately find $\frac{i}{7}$ for every
$i\in\{2,\dots,6\}$. If $k(i)$ is the element in the period that corresponds 
to the digit $i$ in the key, then $\frac{i}{7}=0\boldsymbol{.}k(i)\cdots$ 
where the missing 5 digits are placed as to get a rotation of the 
original period.  In other words, we have 
\begin{equation*}
 \tfrac27= 0.\ol{285714}, \; 
 \tfrac37= 0.\ol{428571}, \; \tfrac47= 0.\ol{571428}, \;
 \tfrac57= 0.\ol{714285}, \; \tfrac67= 0.\ol{857142}.
\end{equation*}
Equivalently, each number $i\big(\frac{10^6-1}{7}\big)$ is a rotation of
$\frac{10^6-1}{7}=142857$. For this reason, $142857$ is called a
\emph{cyclic number}. This property is preserved in some other arithmetic
bases. For instance, in base $3$ we have
\begin{gather*} 
 \tfrac17 = (0.\ol{010212})_3, \;
 \tfrac27 = (0.\ol{021201})_3, \; \tfrac37 = (0.\ol{102120})_3, \\ 
 \tfrac47 = (0.\ol{120102})_3, \;
 \tfrac57 = (0.\ol{201021})_3, \; \tfrac67 = (0.\ol{212010})_3.  
\end{gather*}
Therefore, we say that $7$ generates the $3$-cyclic number $(010212)_3$. 
Similarly, in base $17$, we get  
\begin{gather*} 
 \tfrac17 = (0.\ol{274e9c})_{17}, \;
 \tfrac27 = (0.\ol{4e9c27})_{17}, \; \tfrac37 = (0.\ol{74e9c2})_{17}, \\ 
 \tfrac47 = (0.\ol{9c274e})_{17}, \;
 \tfrac57 = (0.\ol{c274e9})_{17}, \; \tfrac67 = (0.\ol{e9c274})_{17},  
\end{gather*}
where $c=12$ and $e=14$. Thus $7$ also generates a $17$-cyclic number.
However, if we use the binary system, we get 
\begin{gather*}
 \tfrac17 = (0.\ol{001})_{2}, \; \tfrac27 = (0.\ol{010})_{2}, \; 
 \tfrac37 = (0.\ol{011})_{2}, \\
 \tfrac47 = (0.\ol{100})_{2}, \; \tfrac57 = (0.\ol{101})_{2}, \; 
 \tfrac67 = (0.\ol{110})_{2},  
\end{gather*}
so $7$ does not generate a $2$-cyclic number. 

If we are flexible about the base $b$ chosen to represent a number, then 
for every odd prime $p$ we can always pick a $b$ such that $\frac{1}{p}$ 
has a maximal period in that base.
In other words, every prime generates a $b$-cyclic number for some base
$b$. How many such bases are there?  Is there an easy way to change from 
one such base to another?  Are there $b$-cyclic numbers for every $b$?  
Observe that: 
\begin{quote} \em
An odd prime $p$ generates a $b$-cyclic number if and only if $b$ is 
a primitive root mod $p$. 
\end{quote}
Therefore, there are infinitely many bases for which the given prime $p$ 
gives rise to a cyclic number. On the other hand, whether a given $b$ 
(not a square number) is a primitive root for infinitely many primes is 
still an open problem (Artin's conjecture). In any case, in the decimal 
system ($b=10$), the numbers 7, 17, 19, 23, 29, etc. are long primes and 
generate cyclic numbers. However, in the hexadecimal system ($b=16$) 
there are no cyclic numbers at all because a square number is never a 
primitive root.

Now we know that the reason for $7$ to generate cyclic numbers in 
the bases $3$, $10$, and $17$ is because they are primitive roots mod $7$.  
Moreover, we can easily write $\frac17$ in any base of the form $b=3+7t$ 
by using the key $\langle 132645\rangle$.  For instance, the period of 
$\frac17$ in base $10$ can be obtained from the period in base $3$ by the 
rule:
\begin{displaymath} 
\begin{array}{rl}
010212 & \leftarrow \text{base 3}\\
\underline{+\,\it 132645} & \leftarrow \text{key} \\ 
142857 & \leftarrow \text{base 10} 
\end{array}
\end{displaymath}
In the same way,  
\begin{displaymath} 
\begin{array}{rl}
142857 & \leftarrow \text{base 10}\\
\underline{+\,\it 132645} & \leftarrow \text{key} \\ 
274e9c& \leftarrow \text{base 17} 
\end{array}
\end{displaymath}
In fact, the period of $\frac17$ in base $3+7t$ consists of the 6 digits
obtained by adding $t$ times $\langle 132645\rangle$ to the period 
$010212$. In particular, since $2005=3+7\times 286$,
\[ \tfrac17=(0.\overline{[286][859][572][1718][1145][1432]})_{2005}. \]  
This follows from Theorem~\ref{ChangeBase} which is a
consequence of the following lemma.

\begin{lemma}\label{Lem1}
Let $b,m>1$ be integers such that $(b,m)=1$.  Then, in base $b$, we get 
the representation 
$\frac{1}{m}=(0\boldsymbol{.}\, a_1 a_2\cdots a_{i}\cdots )_b$ with
\[ a_i=\frac{b}{m}(b^{i-1}\!\!\mod m) - \frac{1}{m}(b^i\!\! \mod m). \]
\end{lemma}
\begin{proof}
For any $i$, we may write 
\[ a_1\cdots a_i\boldsymbol{.}\,a_{i+1}\cdots=\frac{b^i}{m} = 
   \frac{1}{m}(b^i\!\!\mod m) + \frac{1}{m}(b^i - (b^i\!\!\mod m)). \]
Since $(b,m)=1$, we have $0< \frac{1}{m}(b^i\!\!\mod m)<1$, so the
second term in the latter sum must be the integer part of $\frac{b^i}{m}$. 
Therefore, 
\begin{align*}
a_1\cdots a_{i-1}a_{i} &=\frac{1}{m}(b^i - (b^i\!\!\mod m)), \\
a_1\cdots a_{i-1}0 &=\frac{b}{m}(b^{i-1} - (b^{i-1}\!\!\mod m)).
\end{align*}
The formula for $a_i$ follows now by taking the difference.
\end{proof}

There is a simple algorithm to change the representation of the fraction 
$\frac{1}{m}$ in base $b$ to a base of the form $b+mt$.

\begin{theorem}\label{ChangeBase}
Let $b,m>1$ be integers such that $(b,m)=1$. Assume that $\frac{1}{m}$
is represented in base $b$ as 
$\frac{1}{m}=(0\boldsymbol{.}\, a_1 a_2\cdots a_{i}\cdots )_b$.
Then for any $t\in\mathbb{N}$, the fraction $\frac{1}{m}$ can be
represented in base $(b+mt)$ as 
\[ \tfrac{1}{m}=(0\boldsymbol{.}\, a'_1 a'_2\cdots a'_{i}\cdots)_{b+mt}, \] 
where $a'_i = a_i + tk_i$ with $k_i = (b^{i-1}\!\!\mod m)$.
\end{theorem}
\begin{proof}
By Lemma~\ref{Lem1} we know that for $i\in\mathbb{N}$,
\begin{align*}
a_i &=\frac{b}{m}(b^{i-1}\!\!\mod m) - \frac{1}{m}(b^i\!\! \mod m),\\
a'_i &=\frac{b+mt}{m}((b+mt)^{i-1}\!\!\mod m) - 
       \frac{1}{m}((b+mt)^i\!\! \mod m).
\end{align*}
On the other hand, $(b+mt)^{i-1}\equiv b^{i-1}\!\!\mod m$ and 
$(b+mt)^{i}\equiv b^{i}\!\mod m$. Thus
\begin{equation*}
 a'_i-a_i = t (b^{i-1}\!\!\mod m)=tk_i.
\end{equation*}
\end{proof}

The case $m=7$ in the bases $3$, $10$, and $17$, was discussed above.
A closer look to the proofs of Lemma~\ref{Lem1} and
Theorem~\ref{ChangeBase} reveals that similar statements hold 
even if $(b,m)\not=1$.  In order to illustrate such a situation we will
consider $\frac14$. According to the algorithm described in
Theorem~\ref{ChangeBase}, if we find the representation of $\frac14$ in 
the bases $2$, $3$, $4$, and $5$, together with the corresponding
``keys'', then we will be able to easily represent $\frac14$ in any base.
The key $\langle k_1 \cdots k_\ell\rangle$ associated with $\frac{1}{m}$
in base $b$ is defined by $k_i=(b^{i-1}\!\!\mod m)$ where $\ell$ is either 
the length of the fundamental period of $\frac{1}{m}$ or the length of its 
nontrivial fractional part. Thus
\begin{alignat*}{2}
\tfrac14 &= (0.01)_2 \to \langle 12\rangle & \qquad
\tfrac14 &= (0.\overline{02})_3 \to \langle 13\rangle \\
\tfrac14 &= (0.1)_4 \to \langle 1\rangle & \qquad
\tfrac14 &= (0.\overline{1})_5 \to \langle 1\rangle 
\end{alignat*}
Hence $\frac14=(0.13)_6$ since $01+{\it 12}=13$, and 
$\frac14=0.25$ because $13+{\it 12}=25$. Similarly, using 
$\langle 1\rangle$ as key, one gets for instance
\[ \tfrac14=(0.2)_8, \quad \tfrac14=(0.4)_{16}, 
   \;\text{ and }\; \tfrac14=(0.\overline{3})_{13}. \]
In particular, in base $2005$, $\frac14=(0.\overline{[501]})_{2005}$.

\section{Extension of Midy's theorem}
In the second part of this note we give an extension of Midy's Theorem
(cf. \cite{Dickson}) which in a simplified version states:
\begin{quote}\em
If $p$ is a prime number and the period for $\frac{1}{p}$ has even length, 
then the sum of its two halves is a string of $9$'s. 
\end{quote}
For example, 
\begin{align*}
 \tfrac17 &=0.\overline{142857} \text{ and } 142+857=999,\\
 \tfrac{1}{17} &=0.\overline{0588235294117647} \text{ and } 
 05882352+94117647=99999999.
\end{align*}
It turns out that if the length of the period is divisible by $3$, the
same is true for the sum of its thirds. For instance, in the case of 
$\frac17$, we have indeed $14+28+57=99$. A proof of this 
generalization can be found in \cite{Gin04}.

At this point, two natural questions arise: 
\begin{itemize}
\item What happens in other arithmetic bases? 
\item What about other partitions of the period? 
\end{itemize}
Our intention in this section is to answer these questions, but first 
we need the following lemma that will be crucial to prove 
Theorem~\ref{GeneralMidy}.

\begin{lemma} \label{Lem2}
Let $p$ be an odd prime and let $b>1$ be an integer such that $(b,p)=1$. 
Let $\frac{1}{p}$ be represented in base $b$ as 
$\frac{1}{p}=(0\boldsymbol{.}\,a_1 a_2\cdots a_i\cdots)_b$. 
If $\ell$ is a composite number and $d$ is a divisor, then for every 
$i\in\mathbb{N}$, we have
\begin{equation}\label{digits}
  a_{(i-1)\frac{\ell}{d}+1}\dots a_{i\frac{\ell}{d}}=
  \frac{b^\frac{\ell}{d}}{p}(b^{(i-1)\frac{\ell}{d}}\!\!\mod p) - 
  \frac{1}{p}(b^{i\frac{\ell}{d}}\!\!\mod p),
\end{equation}
where $a_j\cdots a_k$ denotes the integer with digits $a_j$,\ldots, $a_k$.
\end{lemma}
\begin{proof}
As in the proof of Lemma~\ref{Lem1} we have
\[ a_1\cdots a_{j} =\frac{1}{p}(b^j - (b^j\!\!\mod p)) \]
for every $j$. In particular,
\begin{align*}
a_1\cdots a_{(i-1)\frac{\ell}{d}}\cdots a_{i\frac{\ell}{d}} 
&=\frac{1}{p}(b^{i\frac{\ell}{d}} - (b^{i\frac{\ell}{d}}\!\!\mod p)), \\
  a_1\cdots a_{(i-1)\frac{\ell}{d}}
  \underbrace{0\cdots 0}_{\frac{\ell}{d}\text{ digits}} 
&=\frac{b^{\frac{\ell}{d}}}{p}(b^{(i-1)\frac{\ell}{d}} - 
  (b^{(i-1)\frac{\ell}{d}}\!\!\mod p)) 
\end{align*}
and \eqref{digits} follows by taking the difference.
\end{proof}

Assume $\frac{1}{p}=(0\boldsymbol{.}\,\overline{a_1 \cdots 
a_{\ell}})_b$ and pick $\ell$ to be the order of $b$ mod $p$, i.e.,
$\ell$ is the smallest positive integer for which $b^\ell\equiv 1\!\!\mod p$.
In other words, $\ell$ is the length of the fundamental period for 
$\frac{1}{p}$ in base $b$. If $d$ divides $\ell$, then we can write 
\begin{equation}\label{partition}
 \tfrac{1}{p}=(0\boldsymbol{.}\,\overline{A_1\cdots A_d})_b
\end{equation}
with $A_i=a_{(i-1)\frac{\ell}{d}+1}\cdots a_{i\frac{\ell}{d}}$.

\medskip
The main result of this section is the following generalization of Midy's
Theorem.

\begin{theorem}\label{GeneralMidy}
Let $p$ be an odd prime and let $b>1$ be an integer such that $(b,p)=1$. 
Let $\ell$ be the length of the fundamental period for $\frac{1}{p}$ in 
base $b$.  Let $d$ be a divisor of $\ell$ and write $\frac{1}{p}$ as in 
\eqref{partition}. If $d$ is even, working in base $b$ we have
\[ \sum_{i=1}^d A_i = \frac{d}{2}[(b-1)\cdots(b-1)]. \]
In words, the sum of the $d$ parts $A_1,\dots,A_d$ of the period for
$\frac{1}{p}$ in base $b$ is $d/2$ times a number whose digits 
are all $(b-1)$.

Moreover, if $d=3$, then
\[ \sum_{i=1}^d A_i = [(b-1)\cdots(b-1)]. \]
\end{theorem}

\begin{proof}
Let $d$ be a divisor of $\ell$. We use Lemma~\ref{Lem2} to write
\begin{align*}
 \sum_{i=1}^d A_i &= \sum_{i=1}^d 
 \left[\frac{b^\frac{\ell}{d}}{p}(b^{(i-1) \frac{\ell}{d}}\!\!\mod p)
 -\frac{1}{p}(b^{i\frac{\ell}{d}}\!\!\mod p)\right] \\
 &= \frac{b^\frac{\ell}{d}}{p} 
    \sum_{i=1}^d (b^{(i-1) \frac{\ell}{d}}\!\!\mod p) 
   -\frac{1}{p} \sum_{i=1}^d (b^{i\frac{\ell}{d}}\!\!\mod p)\\
 &= \left( \frac{b^\frac{\ell}{d}-1}{p}\right)
    \sum_{i=1}^d (b^{(i-1) \frac{\ell}{d}}\!\!\mod p)
\end{align*}
since $b^\ell\equiv 1\!\!\mod p$.

Note that $b^\frac{\ell}{d}-1 = (b-1)\sum_{j=0}^{\frac{\ell}{d}-1} 
b^j$ is exactly the number $[(b-1)\cdots(b-1)]$ in base $b$. It has
$\ell/d$ digits. Also, for every $1\le i\le d$, $b^{(i-1)\frac{\ell}{d}}$ is a 
unique $d$th root of unity mod $p$, that is, it is a solution of 
$x^d-1\equiv 0\!\!\mod p$. 

If $d$ is even and $r$ is a $d$th root of unity mod $p$, then so is $p-r$. 
Therefore,  
\[ \frac{1}{p}\sum_{i=1}^d (b^{(i-1)\frac{\ell}{d}}\!\!\mod p)
   =\frac{d}{2} \]
and we are done with the first statement of the theorem. 

If $d=3$, then
\[ \frac{1}{p}\sum_{i=1}^3 (b^{(i-1)\frac{\ell}{3}}\!\!\mod p)
   =\frac{1+r_1+r_2}{p}, \]
where $0<r_1,r_2\le p-1$ are the other two third roots of unity mod $p$.
Since $0<1+r_1+r_2\le 1+2(p-1) < 2p$, and since $p$ clearly divides  
$1+r_1+r_2$, we must have $1+r_1+r_2=p$.  So the second statement is proven.
\end{proof} 

\begin{example}
Consider the fraction $\frac{1}{13}$ in base $6$:
\[ \tfrac{1}{13}=(0.\overline{024340531215})_6 \]
If we split the period in 2, 4, or 6 parts, we get 
\begin{gather*}
024340+531215 = 555555 \\
024+340+531+215 = (024+531)+(340+215) = 2\times 555 \\
02+43+40+53+12+15 = (02+53)+(43+12)+(40+15) = 3\times 55.
\end{gather*}
If we split the period in 3 parts, then
\[ 0243+4053+1215 = 5555 \]
noting that the sum must be performed in base $6$.
\end{example}

The general case when the divisor $d$ is odd $>3$ does not seem to obey 
any reasonable pattern. Clearly, we always have
\[ \sum_{i=1}^d A_i = \alpha_d [(b-1)\cdots(b-1)] \]
for some $\alpha_d\ge 1$.  In fact, as shown in our proof, $\alpha_d$ is
$1/p$ times the sum of the $d$th roots of unity mod $p$. Surprisingly,
$\alpha_d$ does not depend on the base $b$.

We finish this note illustrating some examples for $d=5$. We consider
\begin{align*}
\tfrac{1}{11}&=(0.\overline{0001011101})_2 \\	
\tfrac{1}{31}&=(0.\overline{000212111221020222010111001202})_3 \\
\tfrac{1}{101}&=(0.\overline{00504337031261331765\cdots})_8.
\end{align*}
The interesting feature of these examples is that they give different
values for the constant $\alpha_5$.
For $\frac{1}{31}$ we work in base $3$ and get
\[ 000212+111221+020222+010111+001202=222222. \]  
Now, in base 2, we have for $\frac{1}{11}$,
\[ 00+01+01+11+01 = (00+01+01+01) + 11 = 2\times 11. \]
Finally, we look at $\frac{1}{101}$ in base $8$. Again, we break the digits 
of the period into $5$ groups:
\begin{align*}
&00504337031261331765\\
&67101715235114024215\\
&74145305547727344074\\
&65164460121067606254\\
&26637535620363247223\\
\end{align*}
and add these numbers in base $8$. Their sum is $277777777777777777775$ 
which is equal to $3\times 77777777777777777777$. 

In general, $1\le\alpha_5\le 3$. As we just saw, every value in that
range is possible.  Further work on this subject may include finding sharp 
bounds for $\alpha_d$ when $d$ is odd greater then $3$.

{\em Penn State Altoona, 3000 Ivyside Park, Altoona, PA 16601.}\\
{jgil@psu.edu, mdw8@psu.edu}

\end{document}